\documentclass[review]{elsarticle}

\usepackage{lineno,hyperref,subfigure,epstopdf,amsmath,amssymb,color}

\modulolinenumbers[5]

\newcommand{\ds}{\displaystyle}

\newcommand{\xp}{\hat{x}_k^-}

\newcommand{\xf}{\hat{x}_{k(uc)}}
\newcommand{\xc}{\hat{x}_k}
\newcommand{\zz}{\mathcal{Z}}
\newcommand{\Pp}{P_{k(uc)}^-}
\newcommand{\Pc}{P_k}
\newcommand{\Pf}{P_{k(uc)}}
\newcommand{\PCf}{P^C_{k(uc)}}
\newcommand{\PCc}{P^C_k}
\newcommand{\Ks}{K_k^\prime}
\newcommand{\Cc}{\hat{C}_k}
\newcommand{\Cf}{\hat{C}_{k(uc)}}
\newcommand{\bmat}[1]{\begin{bmatrix}#1\end{bmatrix}}
\newcommand{\comment}[1]{}
\newcommand{\Exp}[1]{\mathbf{E}\left[#1\right]}

\journal{Applied Mathematics and Computation}

\begin{document}

\begin{frontmatter}

\title{ Nonlinear Kalman Filtering for Censored Observations}


\author[JA]{Joseph Arthur}
\address[JA]{Department of Statistics, Stanford University, Stanford, California, 94305 USA}

\author[AA]{Adam Attarian}
\address[AA]{MIT Lincoln Laboratory, Massachusetts Institute of Technology, Lexington, Massachusetts, 02420 USA}

\author[NCSU]{Franz Hamilton\corref{cor1}}
\ead{fwhamilt@ncsu.edu}
\cortext[cor1]{Corresponding author}

\author[NCSU]{Hien Tran}

\address[NCSU]{Department of Mathematics and Center for Quantitative Sciences in Biomedicine, North Carolina State University, Raleigh, North Carolina, 27695 USA}

\begin{abstract}
The use of Kalman filtering, as well as its nonlinear extensions, for the estimation of system variables and parameters has played a pivotal role in many fields of scientific inquiry where observations of the system are restricted to a subset of variables. However in the case of censored observations, where measurements of the system beyond a certain detection point are impossible, the estimation problem is complicated. Without appropriate consideration, censored observations can lead to inaccurate estimates. Motivated by the work of \cite{kalmancensored}, we develop a modified version of the extended Kalman filter to handle the case of censored observations in nonlinear systems. We validate this methodology in a simple oscillator system first, showing its ability to accurately reconstruct state variables and track system parameters when observations are censored. Finally, we utilize the nonlinear censored filter to analyze censored datasets from patients with hepatitis C and human immunodeficiency virus.
\end{abstract}


\end{frontmatter}

\linenumbers

\section{Introduction}
The use of data assimilation for the estimation of unobserved model variables and parameters has become standard practice in modern scientific analysis. Kalman filtering \cite{kalman} and its nonlinear extensions such as the ensemble Kalman filter and extended Kalman filter have gained increasing popularity in application to a variety of problems arising from the physical and biological sciences \cite{enkf7,evensen,rabier,huntsauer,cummings,yoshida,stuart,schiffbook,berry2,ghanim,sauerschiff,ghanim2,hamilton,hamilton2}.

Of recent interest in the field of biomedicine has been the use of ordinary differential equations to model viral infection dynamics such as human immunodeficiency virus (HIV) and hepatitis C virus (HCV) \cite{dahari2007hcv,neumann1998hcv,perelson1996hiv}. Such models can provide insights into disease behavior, treatment, and ultimately improve patient outcomes. Their use for the development of patient specific treatment regimens remains an exciting possibility. However, these models are parameterized by a number of unknown parameters and observation of the system is limited to a noisy subset of the dynamical variables. 

Several methodologies have developed to handle this problem of state and parameter estimation from noisy observations. In particular, the use of Kalman filtering for joint state and parameter estimation has been the topic of several recent papers \cite{voss,sitz2002,matzuka2014,matzuka2012nonlinear}. Unfortunately, this estimation process is further complicated when we consider that the assays used in viral studies for data collection often have a detection limit beyond which accurate measurements are impossible. We refer to these data as \emph{censored}. Measurements within the detection limit are considered \emph{uncensored}. Ignoring the censored data can lead to bias in the estimates \cite{allik1}. As such, a proper framework for handling censored observations is required.

Kalman filtering for censored observations has been the topic of several recent works \cite{kalmancensored,allik1,allik2}. Of particular interest is the method proposed in \cite{kalmancensored}. There the authors derived an auxiliary set of equations for the Kalman filter which provided a modified Kalman gain and covariance update formula to allow for correct inference given censored measurements. The underlying assumption though was that the system of interest is linear. Unfortunately the majority of physical systems and the models representing them are nonlinear, such as those describing the dynamics of HCV and HIV. Our goal in this article is to extend the methodology presented in \cite{kalmancensored} to the case of nonlinear system dynamics. We derive a modified version of the extended Kalman filter allowing for the accurate joint estimation of state variables and parameters in nonlinear systems in the presence of censored data.

We validate our proposed nonlinear censored filter first in a synthetic oscillator system where a detection limit for system observation is imposed. We show the fidelity of filter's state variable and parameter reconstruction even when we have partial observability of the system and several of the data are censored. Additionally, we demonstrate the capability of the filter to track system nonstationarity in the form of a drifting parameter whose dynamics are unknown. Motivated by our success in this synthetic example, we consider the difficult problem of state and parameter estimation for clinical viral data. In particular we examine two datasets from an HCV and HIV clinical study, both of which contain numerous censored data in their respective viral load measurements.

In analyzing these clinical datasets, we follow very closely the work done in \cite{arthur} and \cite{attarian2012} for the HCV and HIV data repsectively. There, the authors provided a detailed model identifiability analysis for these datasets and performed estimation using the expectation maximization algorithm \cite{mclachlan}. Our belief is that the filter should not provide more reliable or accurate estimates than those calculated by expectation maximization, in fact they should be comparable. Therefore we treat the results of \cite{arthur,attarian2012} as ``ground truth" and aim to show that the proposed nonlinear censored filter is able to reproduce similar estimates.  The true utility of the filter is that it provides sequential estimation allowing for the online joint estimation of state variables and parameters and the possibility of tracking parameters whose values drift over time, both of which expectation maximization are unable to do. These capabilities are of particular interest in the field of personalized medicine where researchers may be analyzing clinical data whose measurements span over several months or years and an accurate and timely estimate of the current system state is necessary for appropriate treatment or intervention.

\section{Nonlinear Kalman Filtering with Censored Observations}

We assume the following nonlinear system with continuous-time state dynamics and discrete observations
\begin{eqnarray*}
\dot{x}(t) &=& f(t,x) + w(t)\\
z(t_k) &=& h(x(t_k)) + v_k,
\end{eqnarray*}
where $x$ is an $n$ dimensional state vector and $z$ is an $m$ dimensional observation vector. $w$ and $v$ are Gaussian noise terms with covariances $Q$ and $R$ respectively. The estimation of $Q$ and $R$ is key to the success of any filtering methodology. Here, we perform offline tuning of these error covariance matrices to obtain optimal filter performance.

Due to the system nonlinearity, the standard Kalman filter can not be applied directly. Several nonlinear filters have developed, such as the ensemble Kalman filter (EnKF) and extended Kalman filter (EKF) \cite{simon,law}. Here we focus solely on the EKF, which performs a linearization of the system dynamics at each step of the filter. For a detailed derivation of the algorithm see \cite{davidhiv}.

The EKF is a sequential estimator that consists of a prediction and update step. We solve the following system
\begin{eqnarray*}
\dot{\hat{x}} &=& f(t,\hat{x}) \\
\dot{P} &=& PF^T +FP + Q,
\end{eqnarray*}
with initial conditions $\hat{x}_{k-1}$ and $P_{k-1}$ from $t_{k-1}$ to $t_k$ to compute $\xp$ and $P_k^-$, our prior state and covariance matrix estimates. $F$ is the linearization of the system dynamics, namely $F = \nabla f(\hat{x})$. We form the linearization of the observation operator, $H_k = \nabla h(\xp)$, and then implement the standard Kalman update equations to correct our state and covariance estimates
\begin{align*}
\xc =& \xp + K_k\left[z_k - h(\xp) \right] \\
P_k =& \left[I-K_kH_k\right]P_k^-\\
K_k =& P_k^-H_k^T\left[H_kP_k^-H_k^T+R\right]^{-1}.
\end{align*}

\subsection{Filtering with Censored Data}
In the case of censored data, where the true value of the observation beyond a certain lower or upper detection limit is unknown, the estimation problem is complicated.  Treating these censored observations as uncensored measurements leads to inaccurate estimates. In \cite{kalmancensored}, Gabard\'{o}s and Zufiria addressed this problem of state estimation in the presence of censored data in the Kalman filter framework. The authors derived a new set of equations for the filter which appropriately accounts for censored data during the Kalman update step. In this article we extend these ideas to the nonlinear case, deriving an auxiliary set of update equations for the EKF to accurately handle censored data. The derivation included here follows very closely that in \cite{kalmancensored}, though our assumption throughout is that our system of interest is nonlinear.

We use $U_k$ to denote the vector of all uncensored observations up to time $t_k$.
Similarly, let $C_k$ denote the vector of censored observations, each of which lies in some possibly infinite interval $\zz$. For simplicity, we will write $C_k \in \zz$. The filter proceeds at every step by first estimating the state and error covariance ignoring any censored observations. We denote these naive estimates with $\xf$ and $\Pf$ and use $\xc$ and $\Pc$ to denote the final estimates, which are additionally conditioned on the censored observations lying in $\zz$. To calculate the naive estimates, we use a modified gain term:
\[
K_k = \begin{cases}
0 &
\mbox{if } z_k \mbox{ is censored} \\
P_{k(uc)}^-H_k^T\left[H_kP_{k(uc)}^-H_k^T+R\right]^{-1} &
\mbox{otherwise.}
\end{cases}
\]
Therefore when $z_k$ is censored, we have $\hat{x}_{k(uc)}^- = \xf$ and $\Pp = \Pf$, i.e., the predicted values are equal to the naive estimates.

In the case of a censored observation, we calculate the mean and approximate error covariance for the censored observation conditional on the uncensored data, namely
\begin{align*}
\Cf =& h(\xf)\\
\PCf =& \left[H_k \Pf H_k^T + R \right],
\end{align*}
We also compute
\[
P_{k(uc)}^{Cx} = H_k \Pf,
\]
the covariance between the censored observation and the state. Using multivariate Gaussian calculations (see Appendix), the final state and covariance update equations are defined as
\begin{align}
\xc =& \xf + \Ks\left[\Cc - \Cf \right] \label{xc}\\
P_k =& \Pf - \Ks \left[\PCf - \PCc \right]\left(\Ks\right)^T\label{Pc},
\end{align}
where the new gain term is
\begin{equation}
\Ks = P_{k(uc)}^{xC}\left(P^C_{k(uc)} \right)^{-1}\label{Ks}
\end{equation}
and
\[
P_{k(uc)}^{xC} = \left(P_{k(uc)}^{Cx}\right)^T.
\]

Note that $\Cc$ and $\PCc$ are the mean and covariance of the censored observation given the uncensored observations and conditioned on the censored observation lying in $\zz$. This computation is done using the {\tt tmvtnorm} package in R, which computes the mean and covariance of truncated multivariate normal random variables \cite{Wilhelm:2013aa}. After the first censored observation, \eqref{xc}, \eqref{Pc}, and \eqref{Ks} are used as the state and covariance update equations. Additionally though, we must update $\Cf, \PCf$, and $P_{k(uc)}^{Cx}$ at every step of the filter. This update is carried out in two ways, depending on whether or not $z_k$ is censored.

In the censored case we first update the covariance $P_{k-1(uc)}^{Cx}$ to account for the change in state from $t_{k-1}$ to $t_k$. Momentarily abbreviating this covariance as $D$, we solve the system
\begin{align}
\dot{D} =& DF^T \label{dD}\\
\dot{\hat{x}} =& f(\hat{x}) \label{dxhat}
\end{align}
from $t_{k-1}$ to $t_k$ with initial conditions $D(t_{k-1}) = P_{k-1(uc)}^{Cx}$ and $\hat{x}(t_{k-1}) = \hat{x}_{k-1}$. The result of this computation is that $D(t_k)$ is approximately the covariance between $C_{k-1}$ and $x_k$, conditional on only the uncensored observations (see Appendix for details). We call this covariance $P_{k-1,k(uc)}^{Cx}$ and compute the final updated covariance as
\[
P_{k(uc)}^{Cx} = \bmat{P_{k-1,k(uc)}^{Cx} \\ H_k \Pf}.
\]
Now, we update the naive covariance of the censored observations as
\[
\PCf = \bmat{
P_{k-1(uc)}^C & P_{k-1,k(uc)}^{Cx}H_k^T \\
H_k (P_{k-1,k(uc)}^{Cx})^T & P^z_{k(uc)}},
\]
where the covariance of the new observation is
\[
P^z_{k(uc)} = \left[H_kP_{k(uc)}H_k^T+R\right].
\]
Similarly, updating the naive estimate for the censored observations gives
\[
\Cf = \bmat{\hat{C}_{k-1(uc)} \\ h(\xf) }.
\]

In the case that $z_k$ is not censored, the calculations become slightly more complicated. We first use equations \eqref{dD} and \eqref{dxhat} to compute $P_{k-1,k(uc)}^{Cx-}$, which is equivalent to $P_{k(uc)}^{Cx-}$ since $C_k = C_{k-1}$. This predictive covariance can be updated as
\[
P_{k(uc)}^{Cx} = P_{k(uc)}^{Cx-}\left[I - H_k^TK_k^T \right].
\]
The naive expectation of the censored data vector can be updated according to
\[
\Cf = \hat{C}_{k-1(uc)} + P_{k(uc)}^{Cx-}H_k^T(P^z_{k(uc)})^{-1}
\left[z_k - h(\hat{x}_{k(uc)}^-) \right],
\]
which is analogous to the state update equation in the basic Kalman filter. Similarly, we use the equation
\[
\PCf = P^C_{k-1(uc)}  - P_{k(uc)}^{Cx-}H_k^T (P^z_{k(uc)})^{-1} H_k(P_{k(uc)}^{Cx-})^T
\]
to update the error covariance for the censored observations.

Of course, with an increasing number of censored data the above algorithm can become computationally unwieldy due to the increasing dimension of the covariance matrices. In \cite{kalmancensored} the authors reason that previous censored data can be forgotten over time, allowing for a reduction in the algorithm's computational complexity. In particular the columns of the modified gain term $\Ks$ defined in \eqref{Ks}, where each column corresponds to a censored observation, will naturally decay over time to 0 as more data is processed. Additionally if there are a sufficient number of uncensored observations after a censored measurement, the correlation between the censored observation and the state becomes very small. With these ideas in mind, we can introduce approximations to the state and covariance update by removing past censored observations. This in effect reduces the computational complexity of the algorithm and would allow us to only use subsets of the censored observations for a period of time.


\section{State and Parameter Estimation in Oscillator System}
As an demonstrative example, we consider the estimation problem in the following oscillator system
\begin{eqnarray*}
\dot{x}_1 &=& \alpha x_{2}\\
\dot{x}_2 &=& 4-4x_{1},
\end{eqnarray*}
where $\alpha$ is an arbitrary parameter. Our assumption is that only noisy observations of $x_{1}$ sampled at rate $dt = 0.2$ are available. Observations of $x_1$ though are restricted in that any measurement below a value of 0.8 is censored, implying a censored interval of $\left[-\infty,0.8 \right]$. Given these noisy censored observations, our goal is to estimate $x_{1}$ and $x_{2}$ as well as parameter $\alpha$ using the proposed nonlinear censored filter. The estimation of model parameters with Kalman filtering has received considerable attention. A popular approach is the so-called dual estimation method \cite{hcox,kopp,ericwan,voss} which treats the model parameters $q$ as auxiliary state variables that evolve slowly over time. In this article we assume persistent dynamics, namely $\dot{q} = 0$. Using this approach, we assume $\dot{\alpha} = 0$ and form an augmented state vector consisting of the original state variables $x_1$ and $x_2$ and now $\alpha$, thus allowing for simultaneous state and parameter estimation.

\begin{figure}[!t]
\begin{center}
\includegraphics[width = \columnwidth]{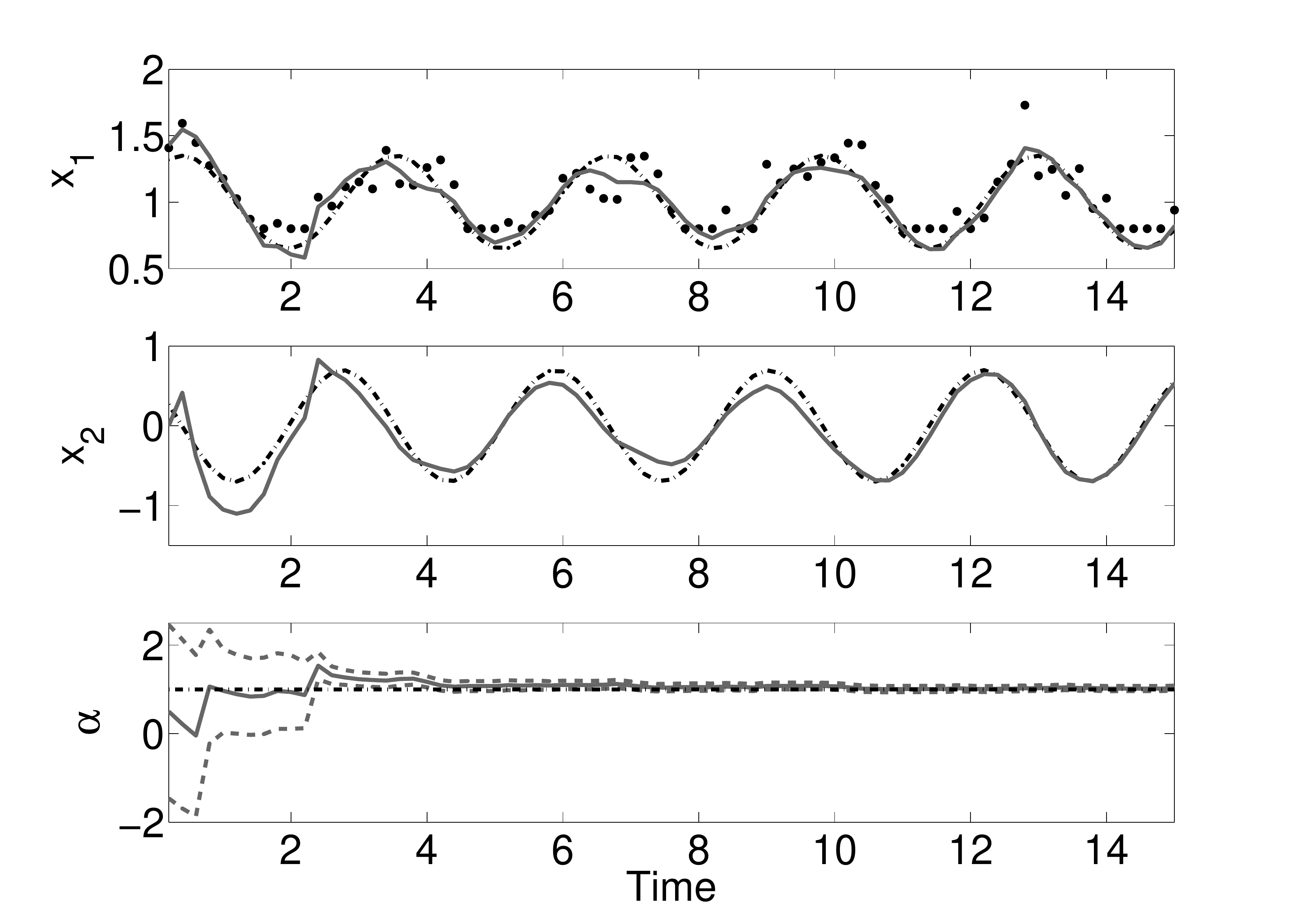}
\end{center}
\caption{State and parameter estimation in oscillator system when $\alpha$ is fixed over time. Observations (black circles) of the $x_{1}$ variable are perturbed by 30\% observational noise and censored below a value of 0.8. Dotted black lines denote the true variable/parameter trajectory and solid grey curves the filter estimates. Dashed grey curve denote the filter estimated 95\% confidence region. Despite the presence of censored data, the filter is able to accurately estimate both state variables as well as the unknown parameter. In particular, we note the fidelity of the reconstruction during censored regions of the data.}
\label{figure1}
\end{figure}

Fig. \ref{figure1} shows the estimation results when $\alpha = 1$ and the observations of $x_{1}$ are corrupted by 30\% Gaussian observational noise. Black circles indicate the noisy observations, dotted black lines denote the true trajectory of the variables and parameters and solid grey curves reflect the filter estimate. In the estimation results for parameter $\alpha$, we also include the filter estimated 95\% confidence interval (dashed grey curves). After an initial transient period, the filter is able to estimate the system variables and parameter with great accuracy. Of particular importance, we notice the fidelity of the filter reconstruction of the variables during the periods of censored data.

\begin{figure}[!t]
\begin{center}
\includegraphics[width = \columnwidth]{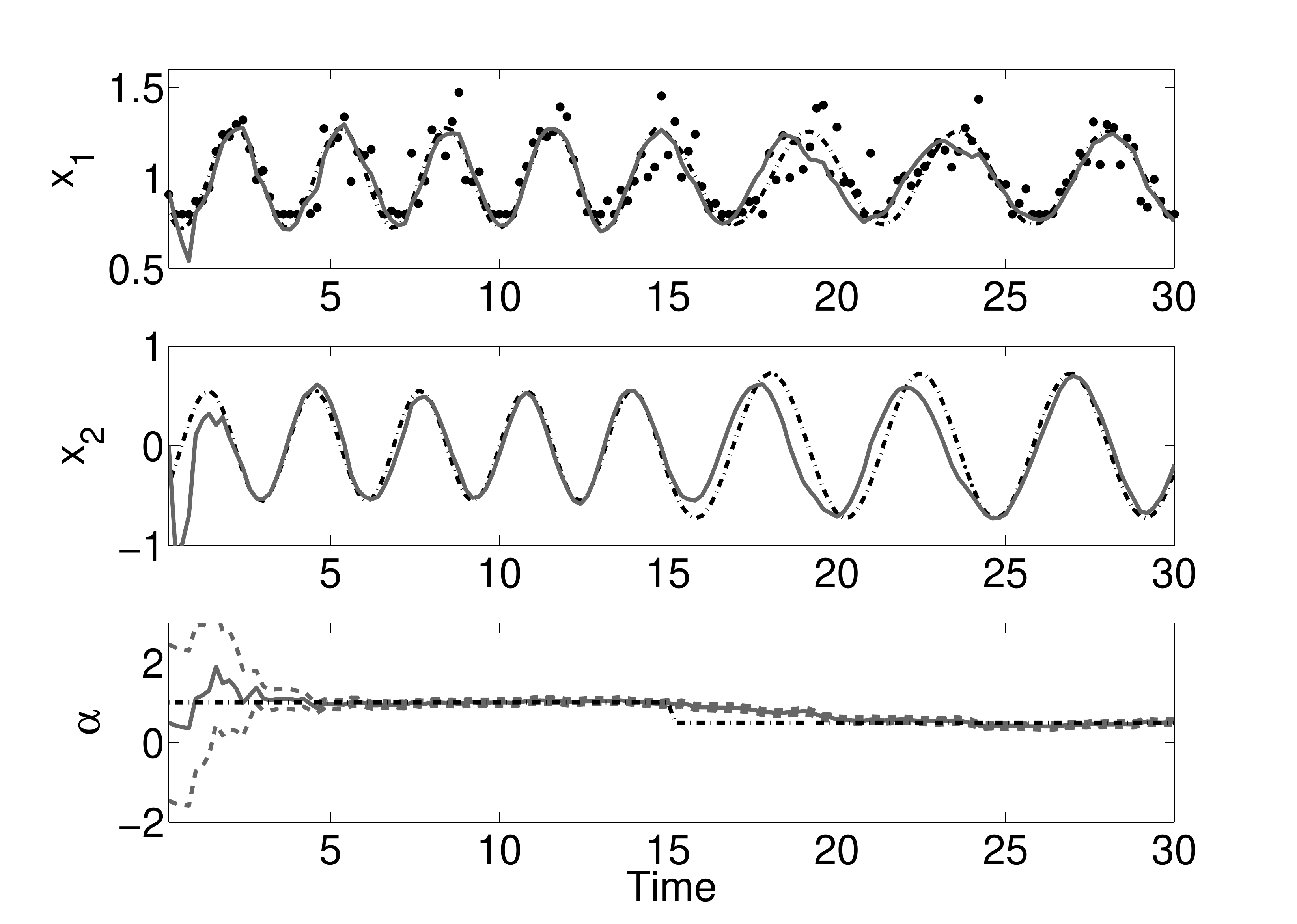}
\end{center}
\caption{State and parameter estimation in simple system when $\alpha$ drifts over time. Observations (black circles) of $x_1$ are perturbed by 30\% Gaussian observational noise. Dotted black lines indicate the true variable/parameter trajectory and solid grey curves the filter estimates. Once again we include the filter estimated 95\% confidence region (dashed grey curve) for $\alpha$. In this more complicated example where a system nonstationarity is present, the filter is still able to accurately track the drift in $\alpha$ and reconstruct the state variables even when there are censored data.}
\label{figure2}
\end{figure}

As previously mentioned, one of the main advantages of using the Kalman filter for estimation is that it is a sequential estimator. While this means that new observations can be processed online without re-analzying the entire dataset, the more important implication is that it allows for the tracking of parameters whose values may drift over time. To simulate this scenario, we considered the estimation problem in the above oscillator system when $\alpha$ changes over time. Namely, its value changes from 1 to 0.5 after 15 units of time. Again, we work under the assumption that only observations of $x_1$ affected by 30\% observational noise are available and also $\dot{\alpha} = 0$. Fig. \ref{figure2} shows the resulting estimation in this nonstationary case. Once again after the initial transient period of the filter we see convergence of $\alpha$ to its correct value and accurate estimation of the state variables. As $\alpha$ drifts, the filter loses track of the $x_1$ and $x_2$ variables but is able to recover after a sufficient amount of data has been observed. Furthermore, the filter is able to accurately track the drift in $\alpha$, despite the presence of censored data.


\section{Estimation in HCV System}
Given the success of the nonlinear censored filter in the oscillator system above, we now consider a significantly more difficult example of state and parameter estimation for analyzing HCV patient data. In this example, an HCV-infected liver undergoes antiviral treatment with interferon-$\alpha$ (IFN) and ribavirin. The typical measurement in this clinical setting is the patient's viral load. Unfortunately viral load is only detectable above a threshold of about 50 copies/mL \cite{snoeck2010hcv}. This means that any measurements below this level are censored (i.e. our censored interval in this case would be $\left[-\infty,50\right]$). Several HCV models have developed, and in particular we consider one by Snoeck, et al. \cite{snoeck2010hcv}. This system is described by the following equations

\begin{eqnarray}
\label{HCVmodel}
\frac{dT}{dt} &=& s+rT\left(1-\frac{T+I}{T_{max}} \right)-dT-\beta V_IT \nonumber \\
\frac{dI}{dt} &=& \beta V_IT+rI\left(1-\frac{T+I}{T_{max}}\right)-\delta I \\
\frac{dV_I}{dt} &=& (1-\bar{\rho} )(1-\bar{\epsilon} )pI - cV_I \nonumber \\
\frac{dV_{NI}}{dt} &=& \bar{\rho} (1-\bar{\epsilon} )pI - cV_{NI} \nonumber,
\end{eqnarray}
where $T$ and $I$ denote concentrations of healthy and infected hepatocytes, and $V_I$ and $V_{NI}$ denote concentrations of infectious and noninfectious virions, respectively. Of note, for parameters $\bar{\rho}$ and $\bar{\epsilon}$ we assume exponentially decaying dynamics
\begin{eqnarray*}
\bar{\rho} &=& \rho e^{-k(t-t_{end})_{+}}\\
\bar{\epsilon} &=& \epsilon e^{-k(t-t_{end})_{+}},
\end{eqnarray*}
where $t_{end}$ indicates the end of treatment and
\begin{gather*}
(a)_+ = 
\begin{cases}
a & \text{if } a\ge0\\
0 & \text{if } a < 0.
\end{cases}
\end{gather*}

%

With respect to (\ref{HCVmodel}), the viral load data maps to the quantity $y = V_I+V_{NI}$. The state variables and parameters of (\ref{HCVmodel}) can range over many orders of magnitude, making accurate estimation difficult. To aid in this process, we apply transformations to all estimated components. In particular, we apply a $\log_{10}$ transformation to the state vector $x$ to compute $\tilde{x}$ and use the relationship
\[
\frac{d\tilde{x}_j}{dt} = \frac{1}{\ln(10)x_j}\frac{dx_j}{dt}
\]
for all $j$, where $x_j = 10^{\tilde{x}_j}$. We also scale the parameters as $\tilde{q}_j = \log_{10}q_j$
for all $q_j$ except the efficacy values $\epsilon$ and $\rho$. These two parameters must be constrained to the interval $[0,1],$ so we instead use
\[
\tilde{q}_j = \tan\left(\pi q_j - \pi/2 \right).
\]
With this log transformation of the model, we assume our filter observation function $h$ to be
\begin{eqnarray*}
h = \log_{10}\left(V_I + V_{NI}\right)
\end{eqnarray*}


The numerous parameters in the model, combined with the limited (and often censored) observability of the physical system, presents a difficult estimation problem. A thorough consideration of this HCV model, including parameter senstivity analysis and estimation, was considered in \cite{arthur}. There, the authors used expectation maximization (EM) to estimate the model states and identifiable parameters for different HCV datasets. Here, we assume the results of \cite{arthur} to be our ``ground truth" and attempt to show that the nonlinear censored filter is able to converge to similar estimates. Again, we emphasize that the censored filter should not give us better or more accurate results than EM, but rather an alternative approach that allows for the sequential estimation of state variables and parameters.

\begin{figure}[!t]
\begin{center}
\includegraphics[width = \columnwidth]{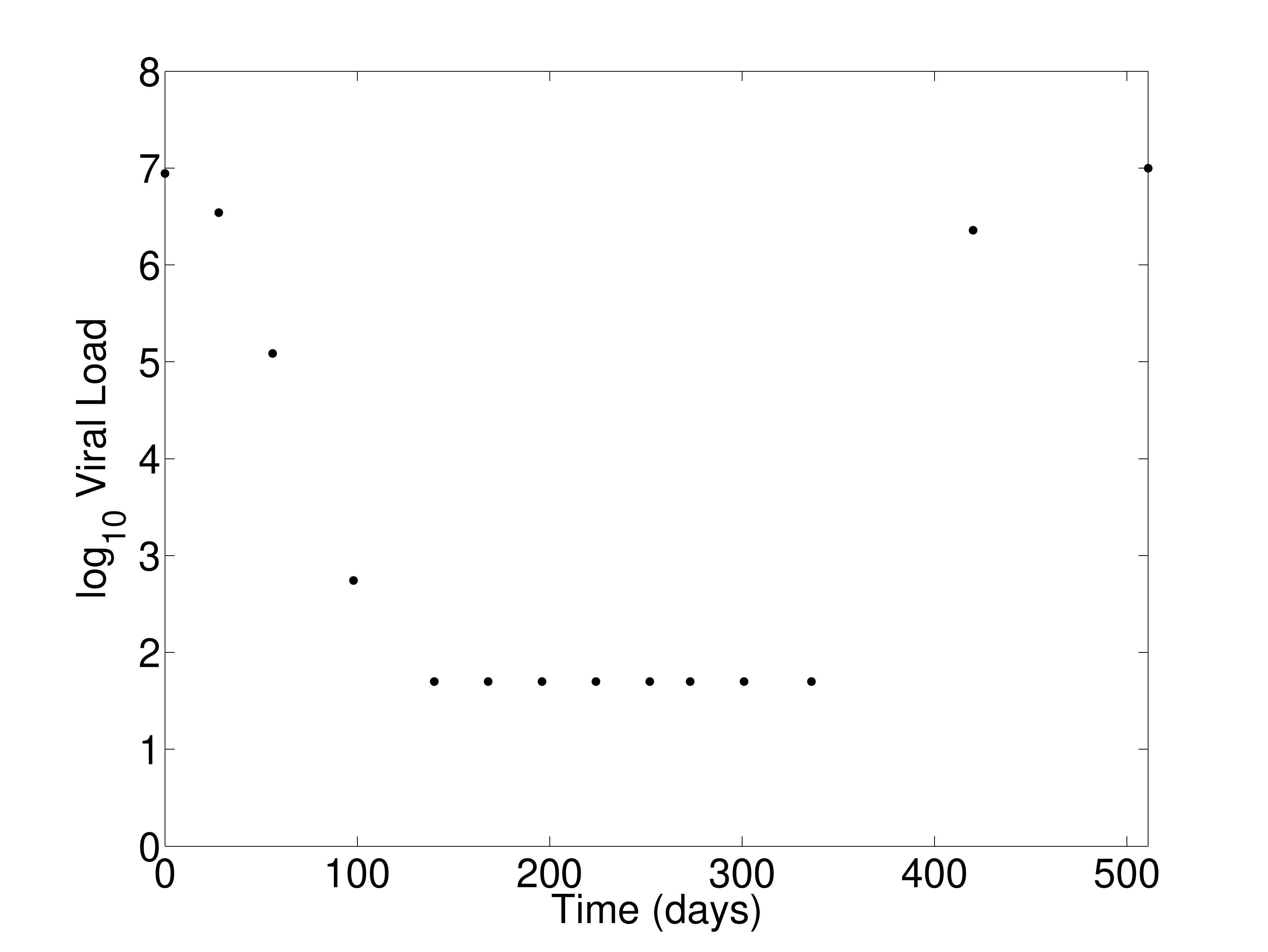}
\end{center}
\caption{Log-scaled viral load data (black circles) collected from a patient in relapse. We observe a clear detection limit in the measurement of viral load, leading to a censored estimation problem.}
\label{figure3}
\end{figure}

For a full description of the model variables and parameters, see \cite{snoeck2010hcv,arthur}. Here, we restrict ourselves to the analysis of data from a patient in relapse as found in \cite{snoeck2010hcv}. Fig. \ref{figure3} shows the log-scaled viral load measurements (black circles) from said patient. We notice immediately that there is a clear lower limit of detection, resulting in censored observability of the system. Following the analysis of \cite{arthur}, we fix the parameter values detailed in Table \ref{tableHCV}. Using the censored filter, we estimate the transformed state variables $\tilde{T}$, $\tilde{I}$, $\tilde{V}_I$ and $\tilde{V}_{NI}$ and parameters $\tilde{\delta}$, $\tilde{c}$ and $\tilde{\epsilon}$ which correspond to the infected cell death rate, virion elimination rate and peginterferon efficacy respectively. Like in our previous example, we assume persistent dynamics for the parameters of interest which allows us to form an augmented state vector and implement the dual-estimation scheme.

\begin{figure}[!t]
\begin{center}
\includegraphics[width = \columnwidth]{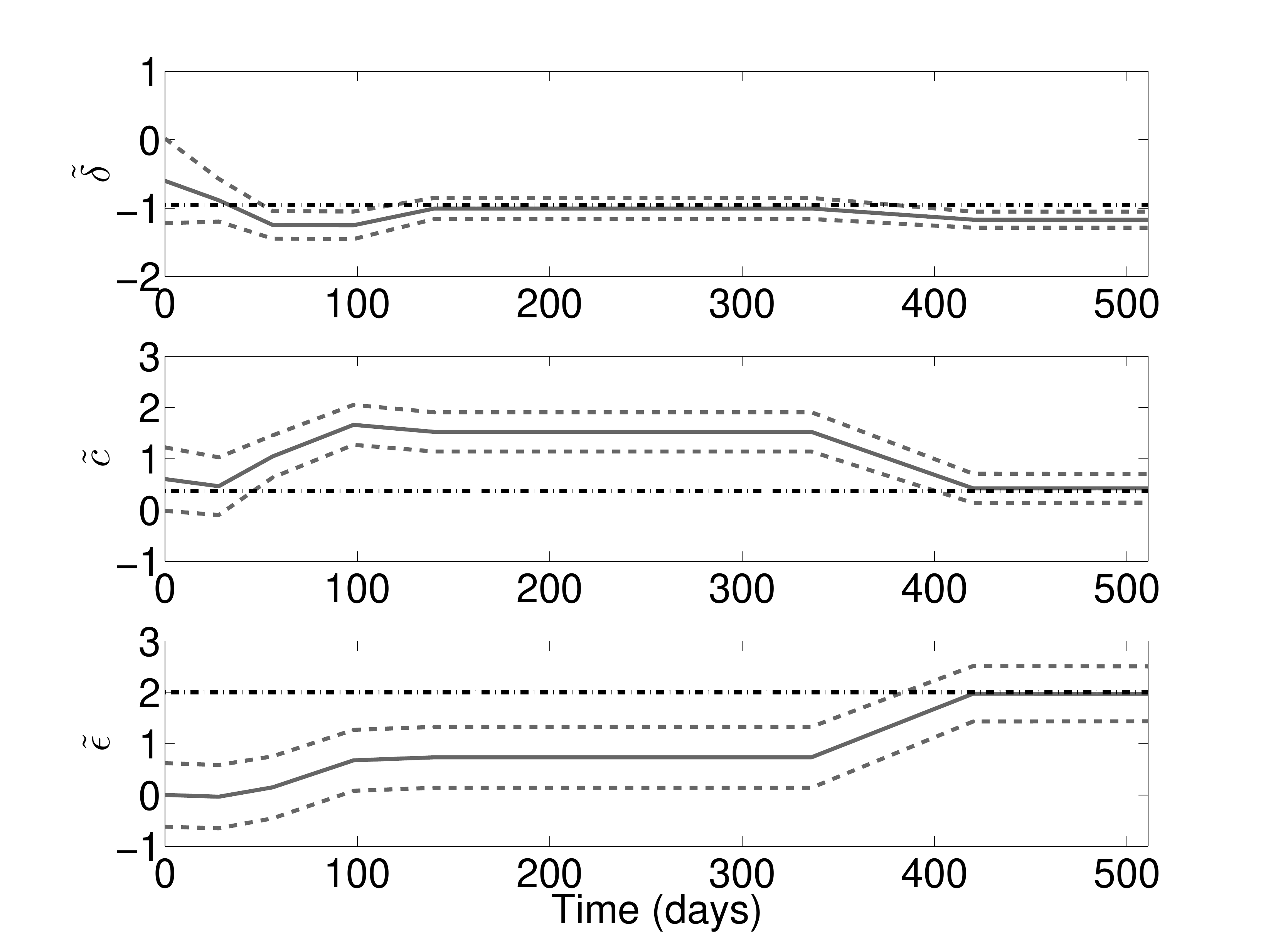}
\end{center}
\caption{Results for the estimation of the transformed HCV model parameters $\tilde{\delta}$, $\tilde{c}$ and $\tilde{\epsilon}$ in relapse data set. Dotted black lines denote ``ground truth" transformed parameter values found in \cite{arthur} and grey solid lines indicate the censored filter estimate. Dotted grey lines indicate the filter estimated 95\% confidence region of the estimate. After an initial transient period, the filter is able to converge to the estimates obtained using EM.}
\label{figure4}
\end{figure}

\begin{figure}[!t]
\begin{center}
\includegraphics[width = 0.9\columnwidth]{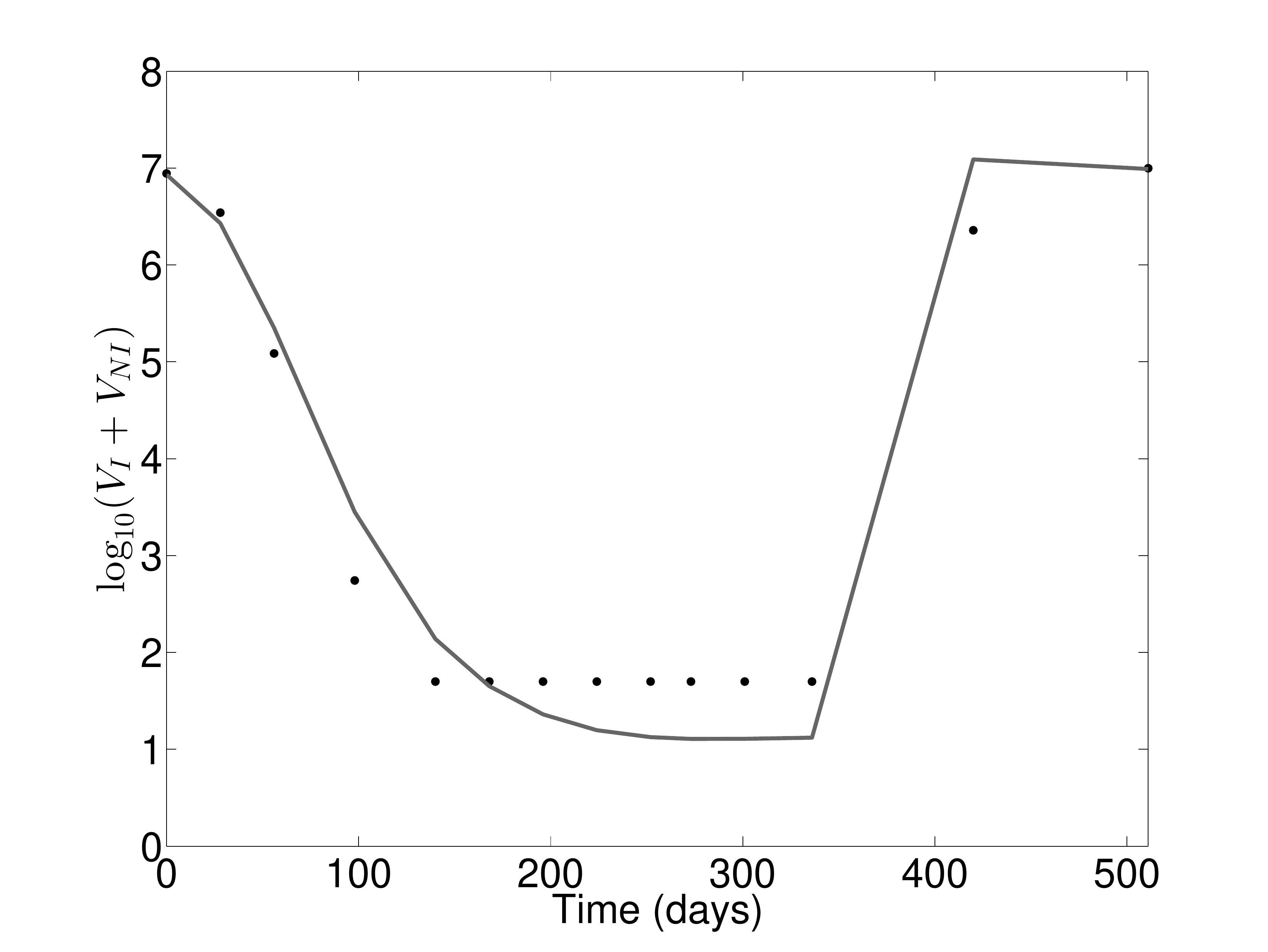}
\end{center}
\caption{Estimation of log-scaled viral load when estimated parameters are fixed to their converged values. Observations (black circles) and filter estimate (grey curve) shown. We obtain a reasonable fit for the data and furthermore estimate a smooth trajectory for the viral load during the censored region of the data.}
\label{figure5}
\end{figure}

Fig. \ref{figure4} shows the results of the HCV parameter estimation for analyzed patient data. Dotted black lines denotes the converged estimate as found using EM in \cite{arthur} and grey curves correspond to the censored filter estimate.  Dotted grey lines indicate the filter estimated 95\% confidence region of the estimates. After a sufficient amount of data, the censored filter is able to converge to parameters estimates comparable to that of EM. Fixing the estimated parameters to their converged estimates, we re-run the filter to obtain an accurate estimation of the state variables. Fig. \ref{figure5} shows the final log-scaled viral load estimate. We obtain a good fit of the data and furthermore we are able to get a reasonable estimation of the system state during the censored data regions.

%

\section{Estimation in HIV System}
We now conisder a more sophisticated example of studying in-host HIV dynamics. The patient data analyzed here comes from a clinical study at Massachusetts General Hospital between 1996 and 2004. This data, originally examined in \cite{adamsthesis}, consists of two measured quantities: CD4$^+$ T-lymphocyte count (cells/$\mu$L) and viral load (copies/mL). Measurement of viral load is once again subject to the detection limits of the assay. In this study both a standard assay, with a detection limit of 400 copies/mL and above, and a high sensitivity assay, with a detection limit of 50 copies/mL and above, were used. Any measurements below the detection limits of the respective assays were effectively censored.


A complex model of in-host HIV dynamics developed in \cite{adams} is described by the following system of equations
\begin{eqnarray}
\label{hivmodel}
\dot{T}_1 &=& \lambda_1-d_1T_1-(1-\epsilon_1(t))d_1V_IT_1 \nonumber \\
\dot{T}_2 &=& \lambda_2-d_2T_2-(1-f\epsilon_1(t))k_2V_IT_2 \nonumber \\
\dot{T}_1^* &=& (1-\epsilon_1(t))k_1V_IT_1-\delta T_1^*-m_1T_1^*E \nonumber \\
\dot{T}_2^* &=& (1-f\epsilon_1(t))k_2V_IT_2-\delta T_2^*-m_2T_2^*E\\
\dot{V}_I &=& (1-\epsilon_2(t))N_T\delta (T_1^*+T_2^*) \nonumber \\
&-&(c+(1-\epsilon_1(t))\rho_1 k_1 T_1 +(1-f\epsilon_1(t))\rho_2 k_2 T_2) V_I \nonumber \\
\dot{V}_{NI} &=& \epsilon_2(t)N_T\delta (T_1^*+T_2^*)-cV_{NI} \nonumber \\
\dot{E} &=& \lambda_{E} + b_E \frac{T_1^*+T_2^*}{T_1^*+T_2^8+K_b}E \nonumber \\
&-&d_E\frac{T_1^*+T_2^*}{T_1^*+T_2^*+K_d}E-\delta_EE \nonumber.
\end{eqnarray}
The model state variables consist of $T_1$ (uninfected type 1 target cells, e.g. CD4$^+$ T-cells), $T_2$ (uninfected type 2 target cells, e.g. magrophages), $T_1^*$ (infected type 1 target cells), $T_2^*$ (infected type 2 target cells), $V_I$ (infectious free virus), $V_{NI}$ (non-infectious free virus) and $E$ (cytotoxic T-lymphocytes, e.g. CD8 cells). The units for the model variables are in $\mu$L. Treatment is modeled through $\epsilon_1(t) = \epsilon_1 u(t)$ and $\epsilon_2(t) = \epsilon_2 u(t)$ where $0 \le u(t) \le 1$.

\begin{figure}[!t]
\begin{center}
\includegraphics[width = \columnwidth]{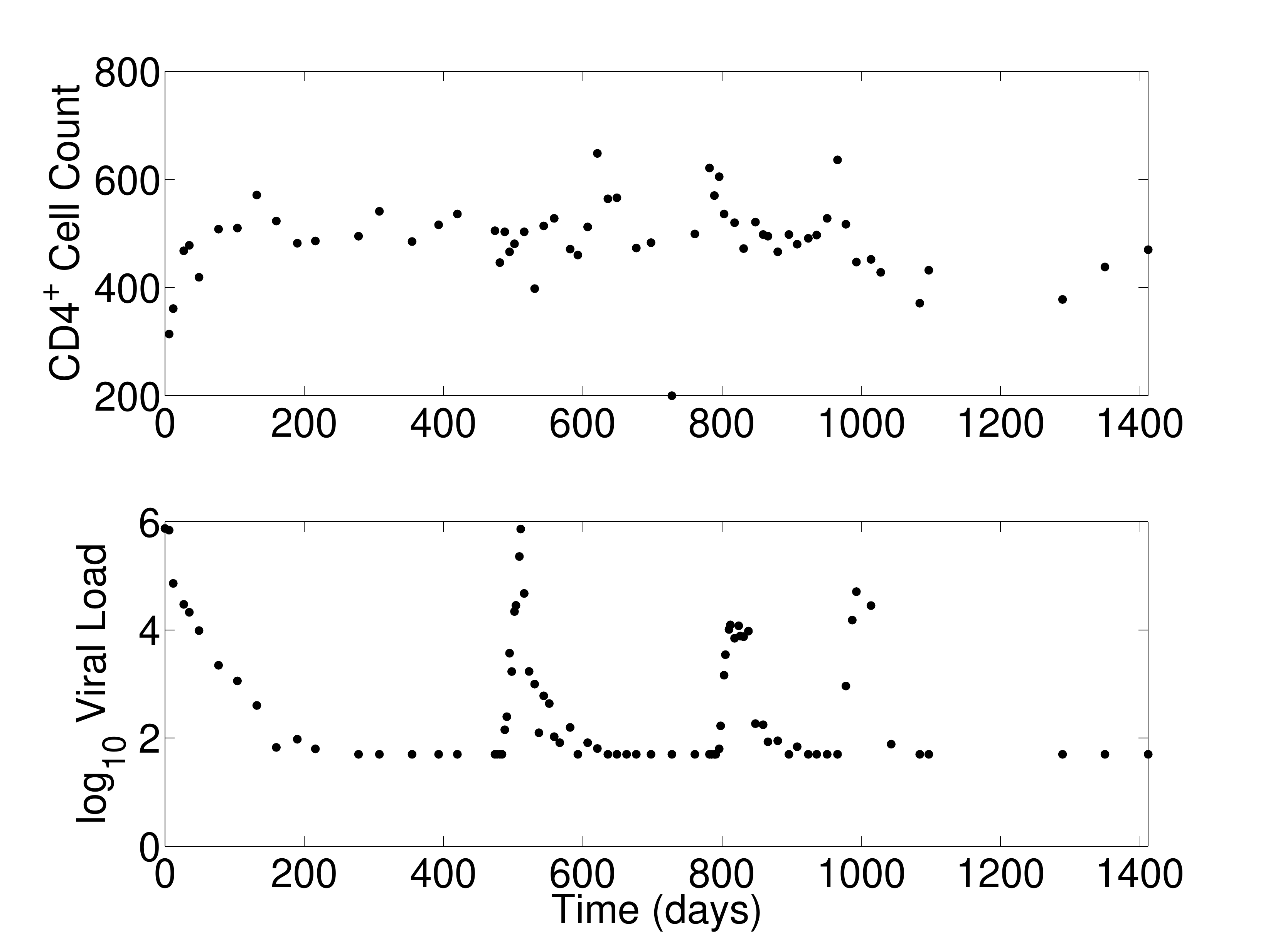}
\end{center}
\caption{Example of patient data (black circles), CD4$^+$ cell count and log-scaled viral load, from clinical study shown. Of note, the measurement of CD4$^+$ and viral load often occur at different time intervals. Furthermore, there is a clear detection limit for measurement of viral load resulting in a censored estimation problem.}
\label{figure6}
\end{figure}

An example of the data collected from a patient in the study is shown in Fig. \ref{figure6}. We notice that the measurement of CD4$^+$ and viral load often occur at different intervals. Additionally, we observe a clear lower limit for viral load detection. With regards to (\ref{hivmodel}), the collected CD4$^+$ data maps to quantity $y_1 = T_1+T_1^*$ and the collected viral load data maps to $y_2 = V_I + V_{NI}$. For a detailed description of (\ref{hivmodel}) and the estimation analysis done for the data acquired in the clinical study, including patient-specific identifiability analysis, see \cite{attarian2012}. Once again, our goal is merely to show that the censored filter derived here is able to reconstruct similar state variable and parameter estimates as those found in \cite{attarian2012} which used the established EM method. We restrict our investigation to the patient data shown in Fig. \ref{figure6}. 

Similarly to the HCV model discussed in the previous section, the HIV model variables and parameters can vary on drastically different orders of magnitude. As such, we once again introduce a $\log_{10}$ transformation for the model variables and parameters to allow for a more robust estimation procedure. The observation function for the filter consists of quantities $h_1 = \log_{10}\left(T_1+T_1^*\right)$ and $h_2 = \log_{10}\left(V_I+V_{NI} \right)$. However as mentioned earlier, the data are collected at different intervals meaning that the filter's observation function changes with respect to the data available at each assimilation time point.

Given these observations, our goal is to estimate log-scaled variables $\tilde{T}_1$,  $\tilde{T}_2$, $\tilde{T}_1^*$, $\tilde{T}_2^*$, $\tilde{V}_I$, $\tilde{V}_{NI}$, $\tilde{E}$ and log-scaled parameters $\tilde{k}_1$ and $\tilde{k}_2$ which correspond to the population 1 and population 2 infection rates respectively. As in our previous examples, we assume persistent dynamics for the parameters thereby allowing us to implement dual estimation. Parameters that were not estimated were fixed to the values in Table \ref{tableHIV} as detailed in \cite{attarian2012}.

\begin{figure}[!t]
\begin{center}
\includegraphics[width = \columnwidth]{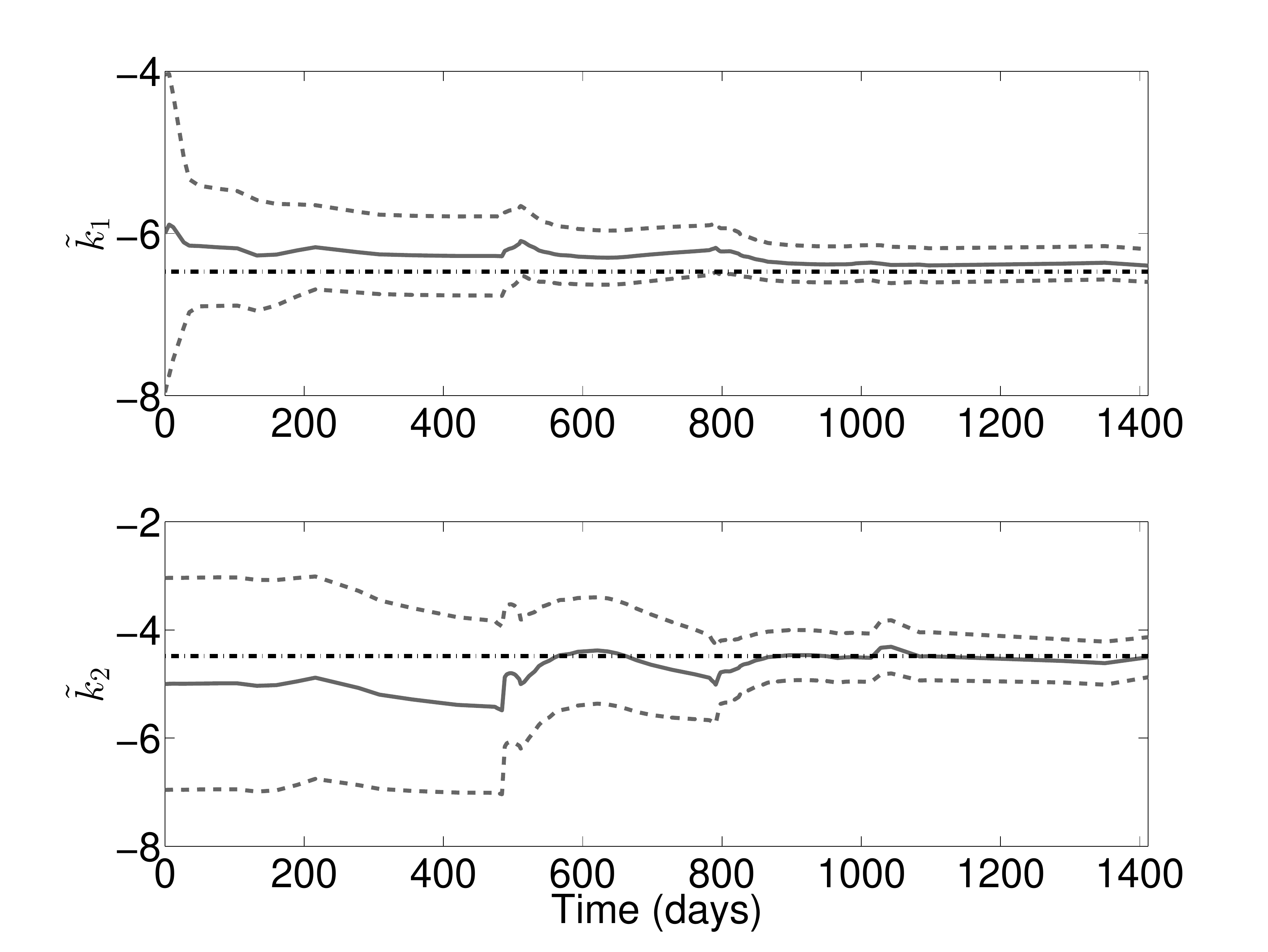}
\end{center}
\caption{Estimated log-transformed parameters for HIV patient dataset. Filter estimated parameter values (grey curve) compare favorably with the values estimated by EM (dotted black line). Filter estimated 95\% confidence interval also shown (dotted grey curves). As the estimates converge to the correct value, the confidence interval shrinks showing reliability of estimates.}
\label{figure7}
\end{figure}

\begin{figure}[!t]
\begin{center}
\includegraphics[width = \columnwidth]{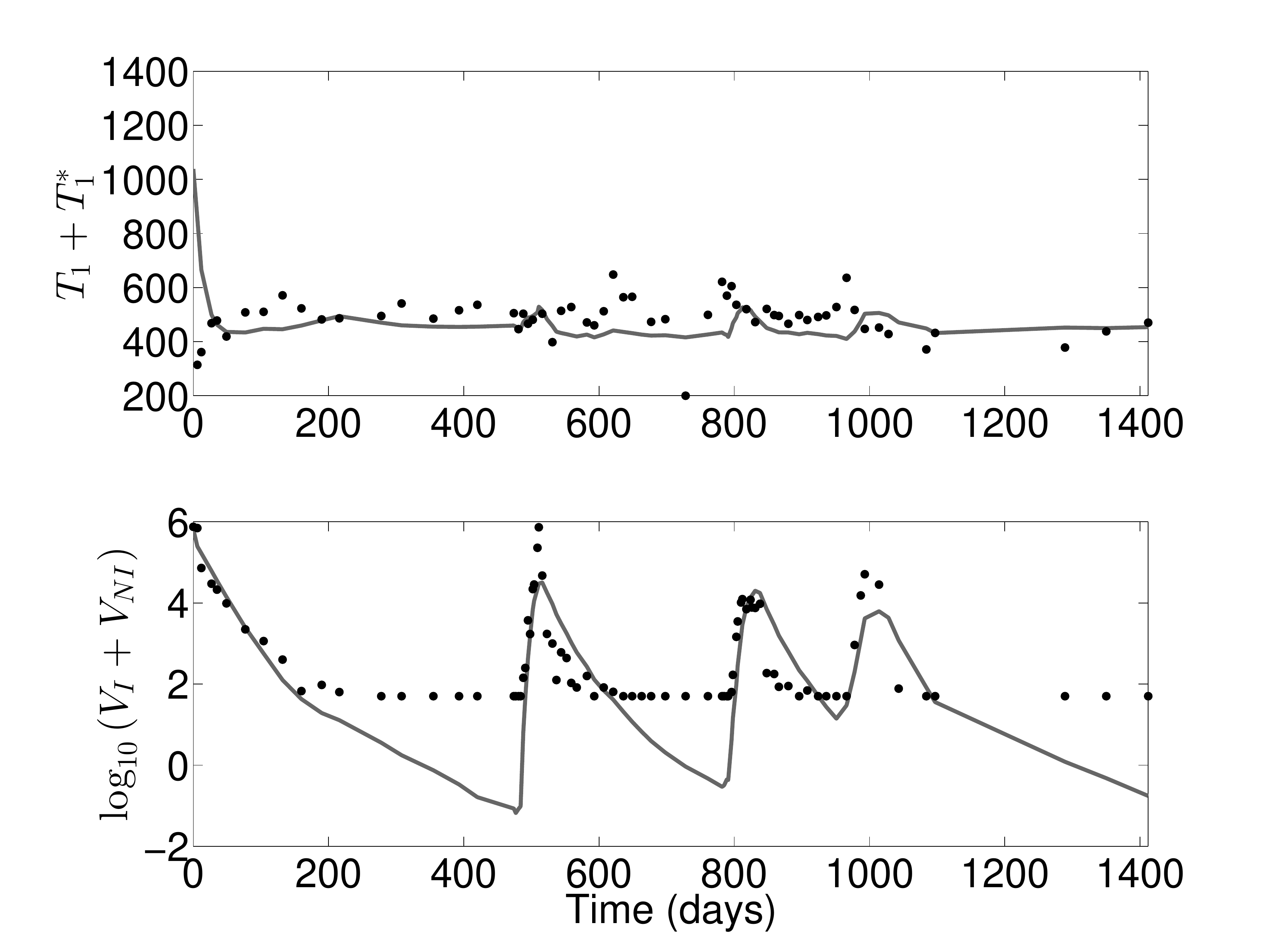}
\end{center}
\caption{Estimation of log-scaled viral load and CD4$^+$ count when estimated parameters are fixed to their converged values. Observations (black circles) and filter estimate (grey curve) shown. We obtain good fits for the data and in particular we obtain a good reconstruction of the viral load during the censored regions of the data.}
\label{figure8}
\end{figure}

Fig. \ref{figure7} shows the results of the filter estimation for transformed parameters $\tilde{k}_1$ and $\tilde{k}_2$. After a sufficient amount of data, the filter estimates (solid grey curves) converge to the desired parameter values (dotted black curve) that were obtained using EM. Additionally, the estimated 95\% confidence region (dashed grey curves) shrinks as convergence occurs reinforcing the optimality of the parameter estimate. We once again fix the estimated parameters to their convergent values and re-run the filter to obtain accurate estimates of the state variables. The resulting filter estimates (solid grey curves) are shown in Fig. \ref{figure8}. We obtain a good fit of both data and additionally are able to get a reasonable reconstruction of the viral load during the regions of censored data, once again showing the capabilities of the censored filter.


\section{Conclusion}
The presence of censored data further complicates the state and parameter estimation process. Incorrectly accounting for these observations can lead to inaccurate estimates resulting in incorrect model inference. Here we derived a modified version of the extended Kalman filter which accounts for the censored observations in the form of an auxiliary set of filter update equations. We examined the performance of this novel filter in an oscillator system where measurements were noisy and censored. We demonstrated its ability to accurately reconstruct state variables and track stationary and drifting parameter values despite the limitations imposed by the censored data. Motivated by this success, we implemented the filter to analyze complex censored data from an HCV and HIV clinical study. The proposed filter was able to obtain comparable estimates for the parameters and state variables as those calculated in the literature using expectation maximization.


The success of the nonlinear censored filter opens up many exciting possibilities. The sequential nature of the algorithm allows for the online estimation of states and parameters and more importantly the tracking of parameters within a patient's dataset that may change over time. Being able to track any potential parameter drift would allow for much more accurate model-based prescription of treatment. Further work should examine the implementation of other nonlinear filters, such as the unscented and ensemble Kalman filter, in place of the EKF which can be costly due to the required system linearization.

\section{Acknowledgments}
This research was partially supported by grants No. RTG/DMS-1246991 and No. DMS-1022688 from the National Science Foundation.

\section{Appendix A. Conditional Moment Calculations}
Suppose $x$ and $z$ are $m$- and $n$-dimensional, jointly Gaussian random vectors. Additionally, let $\zz$ be an $n$-dimensional rectangle in $\mathbb R^n$. Then the conditional mean of $x$ given $z\in\zz$ is
\begin{align*}
\Exp{x|z\in \zz} =& \Exp{\Exp{x|z}|z\in \zz}\\
 =& \Exp{\mu_x + K(z-\mu_z)|z\in \zz}\\
 =& \mu_x + K(\mu_{z|z\in\zz} - \mu_z),
\end{align*}
where $K = P_{xz}P_z^{-1}$ and $\mu_{z|z\in\zz}$ is the conditional mean of $z$ given $z\in\zz$ \cite{kalmancensored}.
The derivation of the conditional covariance is more lengthy, but the result has the simple form
\[
P_{x|z\in\zz} = P_x - K(P_z - P_{z|z\in\zz})K^T,
\]
where $P_{z|z\in\zz}$ is the covariance of $z$ conditional on $z\in\zz$ \cite{kalmancensored}.

\section{Appendix B. Covariance Prediction}
Consider the covariance $D(t)$ between $C$, the vector of censored observations, and the current state vector $x(t)$. This is
\begin{align*}
D(t) =& \Exp{(C-\hat{C})(x(t)-\hat{x}(t))}\\
 =& \Exp{Cx(t)} - \hat{C}\hat{x}(t),
\end{align*}
where, keeping with our censored data Kalman filter, $\hat{C}$ and $\hat{x}(t)$ are expectations given the uncensored data.
We are interested in how $D(t)$ evolves during a time interval when there are no new measurements. Omitting the explicit time-dependence for $x$ and $D$, we have
\begin{align*}
\frac{d}{dt}D =& \frac{d}{dt}\left(\Exp{Cx^T} - \hat{C}\hat{x}^T\right) \\
=& \frac{d}{dt}\Exp{Cx^T} - \hat{C}\frac{d}{dt}\hat{x}^T.
\end{align*}
The first term can be simplified as
\begin{align*}
\frac{d}{dt}\Exp{Cx^T} =& \Exp{C\frac{d}{dt}x^T} \\
=& \Exp{C\left(f(x) + g(t)w(t)\right)^T} \\
=& \Exp{Cf(x)^T} \\
\approx& \Exp{C\left(f(\hat{x}) + \nabla f(\hat{x})(x-\hat{x}) \right)^T}\\
=& \hat{C}f(\hat{x})^T + \Exp{Cx^T}\nabla f(\hat{x})^T - \hat{C}\hat{x}^T\nabla f(\hat{x})^T,
\end{align*}
where we have used the fact that $w(t)$ is independent of $C$ and has expectation $0$. Subtracting off $\hat{C}\frac{d}{dt}\hat{x}^T$ with the substitution $\frac{d}{dt}\hat{x}^T \approx f(\hat{x})^T$ we have
\begin{align*}
\frac{d}{dt}D \approx& \Exp{Cx^T}\nabla f(\hat{x})^T - \hat{C}\hat{x}^T\nabla f(\hat{x})^T\\
=& D\nabla f(\hat{x})^T.
\end{align*}
\comment{
Now, letting $p(x,C)$ be the joint distribution on $x$ and $C$ (conditional on the uncensored data), the first term in the last line is
\begin{align*}
&\ds\int_{-\infty}^\infty\int_{-\infty}^\infty C \left[\frac{d}{dt}x(t)^T\right]p(x,C)dxdC\\
=& \int_{-\infty}^\infty\int_{-\infty}^\infty 
\end{align*}
}
\section*{References}


\begin{table}[!t]
\renewcommand{\arraystretch}{1.3}
\caption{Fixed Parameter Values for HCV Patient Data}
\label{tableHCV}
\begin{tabular}{c|c|c}
Parameter & Description & Value\\
\hline
$\beta$ & Infection rate & $8.7 \times 10^{-9}$\\
$p$ & Virion production rate & 25.1\\
$r$ & Cell proliferation rate & $5.620 \times 10^{-3}$\\
$\rho$ & Ribavirin efficacy & 0.5\\
$k$ & Efficacy decay rate & 0.0238\\
$s$ & Cell production rate & $6.17 \times 10^{4}$\\
$T_{max}$ & Total number of cells per mL & $1.85 \times 10^{7}$\\
$d$ & Cell death rate & 0.003\\
\end{tabular}
\end{table}

\begin{table}[!t]
\renewcommand{\arraystretch}{1.3}
\caption{Fixed Parameter Values for HIV Patient Data}
\label{tableHIV}
\begin{tabular}{c|c|c}
Parameter & Description & Value\\
\hline
$\lambda_1$ & Target cell type 1 source rate & 4.4111\\
$\lambda_2$ & Target cell type 2 source rate & 0.0342\\
$d_1$ & Target cell type 1 death rate & $9.91029 \times 10^{-3}$\\
$d_2$ & Target cell type 2 death rate & $2.6601 \times 10^{-3}$\\
$m_1$ & Population 1 immune-induced clearance rate & $2.8674 \times 10^{-6}$\\
$m_2$ & Population 2 immune-induced clearance rate & $2.9136 \times 10^{-6}$\\
$\rho_1$ & Virions infecting type 1 cell & 0.99052\\
$\rho_2$ & Virions infecting type 2 cell & 0.99622\\
$\delta$ & Infected cell death rate & 0.0952\\
$c$ & Virus death rate & 11.4004\\
$f$ & Treatment efficacy reduction in population 2 & 0.0980\\
$N_T$ & Virions produced per infected cell & 102.5980\\
$\lambda_E$ & Immune effector source rate & $9.4159 \times 10^{-4}$\\
$\delta_E$ & Immune effector death rate & 0.1201\\
$b_E$ & Immune effector max birth rate & 0.0826\\
$d_E$ & Immune effector max death rate & 0.0939\\
$K_b$ & Saturation constant for immune effector birth & 0.1082\\
$K_d$ & Saturation constante for immune effector death & 0.1009\\
$\epsilon_1$ & Reverse transcriptase inhibitor efficacy & 0.5140\\
$\epsilon_2$ & Protease inhibitor efficacy & 0.5770\\
\end{tabular}
\end{table}

\end{document}